\documentclass[11pt]{article}
\usepackage{amscd}
\usepackage{amsthm}
\usepackage{amsmath}
\usepackage{graphicx}
\usepackage{amsfonts}
\usepackage{amssymb}
\newtheorem{theorem}{Theorem}

\newtheorem{corollary}{Corollary}

\newtheorem{definition}{Definition}

\newtheorem{lemma}{Lemma}

\newtheorem{proposition}{Proposition}
\newtheorem{remark}{Remark}

\newtheorem{theorem1}{Theorem\cite{GH}}

\begin{document}

\title{Completable filiform Lie algebras}

\author{Jos\'{e} Mar\'{\i}a Ancochea Berm\'udez\thanks{corresponding author : e-mail : Jose\_Ancochea@mat.ucm.es} and Rutwig
Campoamor\\
Departamento de Geometr\'{\i}a y Topolog\'{\i}a\\
Fac. CC. Matem\'{a}ticas Univ. Complutense\\
28040 Madrid ( Spain )}

\date{}
\maketitle

\begin{abstract}
We determine the solvable complete Lie algebras whose nilradical is isomorphic to a filiform Lie algebra. Moreover we show that for any positive integer $n$ there exists a solvable complete Lie algebras whose second cohomology group with values in the adjoint module has dimension at least $n$.
\end{abstract}

\section{Complete Lie algebras}

Complete Lie algebras first appeared in 1951, in the context of ScheNkman's
theory of subinvariant Lie algebras \cite{Sch}. In his famous derivation tower theorem,
he proved that for centerless solvable Lie algebras the tower of derivation
algebras $Der^{n}\frak{g}=Der^{n-1}\left( Der\frak{g}\right) $ was finite,
and that the last term, which is also centerless, had only inner
derivations. However, the formal definition of complete Lie algeba was not
given until 1962. Their relation with the
tower theorem has proven the importance of these algebras, which also play a
role in rigidity theory of algebraic structures. In the last years, different authors have concentrated on classifications and structural properties of complete Lie algebras \cite{AC3,C,Me1,Me2}, and  various interesting results concerning these structures have been obtained.

\bigskip
For our purpose it will be more convenient to define completeness of Lie algebras in terms of cohomology. Therefore we recall the elementary facts:
 Let $\frak{g}$ be a Lie algebra. A $p$-dimensional
cochain of $\frak{g}$ (with values in $\frak{g}$) is a $p$-linear alternating mapping of $\frak{g}^{p}$ in $\frak{g}$ $\left(  p\in\mathbb{N}^{\ast}\right)
$. A $0$-cochain is a constant function from $\frak{g}$ to $\frak{g}%
$.\newline We denote by $C^{p}\left(  \frak{g},\frak{g}\right)  $ as the space
of the $p$-cochains and
\[
C^{\ast}\left(  \frak{g},\frak{g}\right)  =\oplus_{p\geq0}C^{p}\left(
\frak{g},\frak{g}\right)  .
\]
Over the space $C^{\ast}\left(  \frak{g},\frak{g}\right)  $ we define the
endomorphism
\[
\delta:C^{\ast}\left(  \frak{g},\frak{g}\right)  \longrightarrow C^{\ast
}\left(  \frak{g},\frak{g}\right)
\]%
\[
\Phi\longrightarrow\delta\Phi
\]
by putting
\[
\delta\Phi\left(  X\right)  =X.\Phi\ \ \ \mbox{if}\ \ \ \Phi\in C^{0}\left(
\frak{g},\frak{g}\right)
\]%
\[
\delta\Phi\left(  X_{1},...,X_{p+1}\right)  =\sum_{1\leq s\leq p+1}\left(
-1\right)  ^{s+1}\left(  X_{s}.\Phi\right)  \left(  X_{1},...,\overset
{}{\widehat{X}}_{s},...,X_{p+1}\right)  +
\]%
\[
+\sum_{1\leq s\leq t\leq p+1}\left(  -1\right)  ^{s+t}\Phi\left(  \left[
X_{s},X_{t}\right]  ,X_{1},...,\overset{}{\widehat{X}_{s}},...,\overset
{}{\widehat{X}}_{t,...},X_{p+1}\right)
\]
if $\Phi\in C^{p}\left(  \frak{g},\frak{g}\right)  $, $p\geq1$.\newline By this definition, $\delta\left(  C^{p}\left(  \frak{g}%
,\frak{g}\right)  \right)  \subset C^{p+1}\left(  \frak{g},\frak{g}\right)  $
and we can verify that
\[
\delta\circ\delta=0.
\]
We denote by
\[
\left\{
\begin{array}
[c]{cc}%
Z^{p}\left(  \frak{g},\frak{g}\right)  =Kerd\left|  _{C^{p}\left(
\frak{g},\frak{g}\right)  }\right.   & p\geq1\\
B^{p}\left(  \frak{g},\frak{g}\right)  =Imd\left|  _{C^{p}\left(
\frak{g},\frak{g}\right)  }\right.   & p\geq1
\end{array}
\right.
\]
and $H^{p}\left(  \frak{g},\frak{g}\right)  =Z^{p}\left(  \frak{g}%
,\frak{g}\right)  \;|\;B^{p}\left(  \frak{g},\frak{g}\right)  ,$ $p\geq1$.
 This quotient space is called the cohomology space of $\frak{g}$ of
degree $p$ with values in $\frak{g}$ \cite{Ko}. For $p=0$, then we put $B^{0}\left(
\frak{g},\frak{g}\right)  =\left\{  0\right\}  $ and $H^{0}\left(
\frak{g},\frak{g}\right)  =Z^{0}\left(  \frak{g},\frak{g}\right)  $. This last
space can be identified to the space of all $\frak{g}$-invariant elements that
is
\[
\left\{  X\in\frak{g}\ \mbox{such that}\ adY\left(  X\right)  =0\ \ \forall
Y\in\frak{g}\right\}  .
\]
Then $Z^{0}\left(  \frak{g},\frak{g}\right)  =Z\left(  \frak{g}\right)  $ (the
center of $\frak{g}$).\newline Further, we have
\[
Z^{1}\left(  \frak{g},\frak{g}\right)  =\left\{  f:\frak{g}\longrightarrow
\frak{g\;}|\;\delta f=0\right\}  .
\]
But $\delta f\left(  X,Y\right)  =\left[  f\left(  X\right)  ,Y\right]
+\left[  X,f\left(  Y\right)  \right]  -f\left[  X,Y\right]  $.\ Then
$Z^{1}\left(  \frak{g},\frak{g}\right)  $ is nothing but the algebra of
derivation of $\frak{g}$:
\[
Z^{1}\left(  \frak{g},\frak{g}\right)  =Der\frak{g}.
\]
It is the same for :
\[
B^{1}\left(  \frak{g},\frak{g}\right)  =\left\{  adX,X\in\frak{g}\right\}  .
\]
Thus the space $H^{1}\left(  \frak{g},\frak{g}\right)  $ can be interpreted as
the set of the outer derivations of the Lie algebra $\frak{g}$.

\begin{definition}
A Lie algebra $\frak{g}$ is called complete if $H^{0}\left( \frak{g},\frak{g}%
\right) =H^{1}\left( \frak{g},\frak{g}\right) =\{0\}.$
\end{definition}

Therefore complete Lie algebras are centerless with only inner derivations. Classical examples are the semisimple Lie algebras, and by the Levi decomposition, the class to be analyzed respect to completeness is the one formed by solvable Lie algebras.

\bigskip

Let $\frak{n}$ be a nilpotent Lie algebra and $\frak{t\subset }Der\left( 
\frak{n}\right) $ a maximal toral subalgebra, i.e., an abelian subalgebra
whose elements are $ad$-semisimple and which is maximal for the inclusion
relation. As known, the algebra $\frak{t}$ induces a
root space decomposition over $\frak{n}$ :
\begin{equation*}
\frak{n}=\sum_{\alpha \in \frak{t}^{\ast }}\frak{n}_{\alpha }
\end{equation*}
where $\frak{t}^{\ast }=Hom\left( \frak{t},\mathbb{C}\right) $ and the root
spaces are given by 
\begin{equation*}
\frak{n}_{\alpha }=\left\{ X\in \frak{n\;}|\;[h,X]=\alpha \left( h\right)
.X\;\;\forall h\in \frak{t}\right\} 
\end{equation*}
If $\frak{n}_{\alpha }\neq 0,$ the form $\alpha $ is called a root of $\frak{%
n}$ \cite{Fa}. Let $\Phi =\left\{ \alpha \;\ |\;\frak{n}_{a}\neq 0\right\} $ be the
root system of $\frak{n}$ (associated to the torus $\frak{t}$ ).

\begin{proposition}
Let $\frak{t}\subset Der\left( \frak{n}\right) $ be a torus of $\frak{n}$.
Then there exists a system of generators $\left\{ X_{1},..,X_{n}\right\} $
of $\frak{n}$ and forms $\left\{ \alpha _{1},..,\alpha _{n}\right\} \subset 
\frak{t}^{\ast }$ such that $[h,X_{i}]=\alpha _{i}\left( h\right)
X_{i},\;1\leq i\leq n$.
\end{proposition}

\begin{proof}
The action of $\frak{t}$ over $\frak{n}$ shows that the mapping $\rho :\frak{%
t}\rightarrow Der\left( \frak{n}\right) $ is a representation of $\frak{n}$,
and as the root spaces are irreducible, the representation is completely
reducible. As the commutator algebra $C^{1}\frak{n}$ is clearly a $\frak{t}$%
-submodule of $\frak{n}$, there exists a complementary submodule to $C^{1}%
\frak{n}$. This must generate the algebra, as the generators of a Lie
algebra are extracted from the quotient by the derived algebra. We can
therefore find a basis $\left\{ X_{1},..,X_{n}\right\} $ of the
complementary $\frak{m}$ which satisfies the assumption. The minimality of
the generator system is obvious.
\end{proof}

\begin{definition}
The dimension $r\left(\frak{n}\right)$ of a maximal torus of the preceding form is called the rank of 
$\frak{n}$. 
\end{definition}

The rank is an invariant of the nilpotent Lie algebra. By abuse of notation, we often say that the semidirect product $\frak{n}\oplus\frak{t}$ is of rank $r$, where $\frak{t}$ is a maximal torus of $\frak{n}$ and $dim\left(\frak{t}\right)=r$. A system of generators $\left\{ X_{1},..,X_{n}\right\} $ satisfying the proposition will be called a $t$-msg. 

\begin{theorem}
Let $\frak{g}$ be a solvable complete Lie algebra. Then

\begin{enumerate}
\item  $\frak{g}$ decomposes as $\frak{g}=\frak{t}+\frak{n}$, where $\frak{n}
$ is the nilradical of $\frak{g}$ and $\frak{t}$ is a maximal toral
subalgebra.

\item  the restriction $ad_{\frak{t}}h|_{\frak{n}}$ is an abelian subalgebra
of $Der\left( \frak{n}\right) $, and a maximal torus over $\frak{n}$.
\end{enumerate}
\end{theorem}

This result, whose proof can be found in \cite{Me3}, allows to refer to the
decomposition directly, which is of interest for the structural analysis of
such a Lie algebra, since it shows that we can restrict ourselves to the class of decomposable Lie algebras \cite{AC1}. A nilpotent Lie algebra $\frak{n}$ will be called completable if $\frak{n}\oplus\frak{t}$ is complete, where $\frak{t}$ is a maximal torus of derivations of $\frak{n}$. 
Recall that the minimal number of generators is given by the first Betti number $b_{1}\left(\frak{n}\right)=dim\left(\frac{\frak{n}}{C^{1}\frak{n}}\right)$. 
We in particular obtain an interesting result for maximal rank Lie algebras, i.e., for nilpotent Lie algebras satisfying $dim\left(\frak{t}\right)=b_{1}\left(\frak{n}\right)$:

\begin{theorem}
Let $\frak{n}$ be a nilpotent Lie algebra of maximal rank and $\frak{t}$ a
maximal torus. Then the semidirect product $\frak{t}\oplus\frak{n}$ is a solvable complete Lie algebra.
\end{theorem}

This theorem, due to Zhu and Meng \cite{Me1}, has important consequences for the
theory of complete Lie algebras. We ennounce a couple :

It is easy to see that two solvable complete Lie algebras of maximal rank are isomorphic if and only if their corresponding nilradicals are isomorphic. Moreover, it can be shown that if the nilradical is split, then it is of maximal rank if and only if any component is also of maximal rank. This allows to reduce the study to the determination of complete solvable
algebras whose nilradical is nonsplit. Observe the analogy with the
classical theory.

\begin{remark}
It should be remarked that there exist complete Lie algebras whose
nilradical is not of maximal rank. In low dimensions they have been
classified, but for the general case there only exist partial results
\end{remark}

From the action of the torus on $\frak{n}$ we deduce immediately that central elements must have non-zero eigenvalues, since otherwise we easily construct an outer semi-simple derivation. 

\begin{theorem}
Let $\frak{g}$ be a Lie algebra and $\frak{h}$ a Cartan subalgebra. Assume
that following conditions are satisfied:

\begin{enumerate}
\item  $\frak{h}$ is abelian

\item  $\frak{g=}\frak{h}\oplus\sum_{\alpha \in \Delta }\frak{g}_{\alpha }$ with $\Delta
\subset \frak{h}^{\ast }-\{0\}$

\item[3']  there is a generating system $\left\{ \alpha _{1},..,\alpha
_{l}\right\} $ of $\frak{h}^{\ast }$ in $\Delta $ such that $\dim \frak{g}%
_{\alpha _{j}}=1$ for all $j$ and $\frak{h},\frak{g}_{\alpha _{1}},..,\frak{g%
}_{\alpha _{l}}$ generate $\frak{g}$

\item[4']  Let $0\neq x_{j}\in \frak{g}_{\alpha _{j}}$ and a basis $\left\{
\alpha _{1},..,\alpha _{r}\right\} $ of $\frak{h}^{\ast }$. For $r+1\leq
s\leq l,$%
\begin{equation*}
\alpha _{s}=\sum_{i=1}^{t}k_{is}\alpha _{j_{i}}-\sum_{i=1+t}^{r}k_{is}\alpha
_{j_{i}}
\end{equation*}
where $k_{is}\in \mathbb{N\cup \{}0\},\left( j_{1},..,j_{r}\right) $ is a
permutation  of $\left( 1,..,r\right) $, and there is a formula
\begin{eqnarray*}
&&\lbrack \underset{k_{1s}}{\underbrace{x_{j_{1}},..,x_{j_{1}}}},..\underset{%
k_{ts}}{\underbrace{x_{j_{t}},..,x_{j_{t}}},}...,x_{k_{m}}] \\
&=&[\underset{k_{t+1s}}{\underbrace{x_{j_{t}+1},..,x_{j_{t}+1}}},..\underset{%
k_{rs}}{\underbrace{x_{j_{r}},..,x_{j_{r}}},x_{s},x_{k_{1}}}...,x_{k_{m}}]
\end{eqnarray*}
without regard to the order or the way of bracketing, where $m\neq 0$ if $t=r
$.
\end{enumerate}
Then $\frak{g}$ is a complete Lie algebra.
\end{theorem}

See \cite{Me3} for a proof.

\section{Filiform Lie algebras}
Filiform Lie algebras are the best known and most studied among the nilpotent Lie algebras \cite{AC2,AG2,Ve}. The importance of this class was pointed out 
from 1966 on, when Vergne used these algebras to derive important results on the irreducible components of the variety of nilpotent Lie algebra laws. Various authors have enlarged these results, which explains that this class is the best understoood among the nilpotent ones \cite{C1}. In this section we recall the elementary facts about filiform Lie algebras which will be used.

Let $\frak{n}$ be a nilpotent finite dimensional Lie algebra. Let $%
Y\in \frak{n}-C^{1}(\frak{n})$ be a vector of $\frak{n}$ . Consider
the ordered sequence 
\begin{equation*}
c(Y)=(h_{1},h_{2},\cdots ,)
\end{equation*}
$h_{1}\geq h_{2},...,\geq h_{p},$ where $h_{i}$ is the dimension of the $%
i^{th}$ Jordan block of the nilpotent operator $adY$. As $Y$ is an
eigenvector of $adY$, $h_{p}=1$. Let $Y_{1}$ and $Y_{2}$ be in $\frak{n}-%
D^{1}(\frak{n})$. Let be $c(Y_{1})=(h_{1},...,h_{p_{1}}=1)$ and $%
c(Y_{2})=(k_{1},...,k_{p_{2}}=1)$ the corresponding sequences. We will say
that $c(Y_{1})\geq c(Y_{2})$ if there is an $i$ such that $h_{1}=k_{1},$ $%
h_{2}=k_{2},$ ... , $h_{i-1}=k_{i-1},$ $h_{i}>k_{i}.$ This defines a total
order relation on the set of sequences $c(Y)$.

\begin{definition}
The characteristic sequence of the nilpotent Lie algebra $\frak{n}$ is : 
\begin{equation*}
c(\frak{n})=Sup\{c(Y),Y\in \frak{n}-C^{1}(\frak{n})\}
\end{equation*}
\end{definition}

It is an invariant of $\frak{n}$. A vector $X\in \frak{n}$
such that $c(X)=c(\frak{n})$ is called a characteristic vector of $\frak{n}$.

\begin{definition}
A $n-$dimensional nilpotent Lie algebra $\frak{n}$ whose characteristic
sequence is $c(\frak{n})=(n-1,1)$ is called filiform.
\end{definition}

\begin{lemma} 
If $\frak{n}$ is a filiform Lie algebra, then $r(\frak{n})\leq 2$.
\end{lemma}

For the proof, see \cite{GH}. On the other hand, it is not difficult to prove that filiform Lie algebras satisfy $b_{1}\left(\frak{n}\right)=2$, in view of its characteristic sequence and vector. Therefore it is convenient to separate them by rank.

\begin{proposition}
Every filiform Lie algebras of rank 2 is isomorphic to $L_{n}$ or $%
Q_{n}$, where $L_{n}$ and $Q_{n}$ are the $n$-dimensional Lie algebras defined by 
\begin{equation*}
L_{n}:\left\{ [Y_{1},Y_{j}]=Y_{1+j},\quad j=2,...,n-1\right.
\end{equation*}
\begin{equation*}
Q_{n}=\left\{ 
\begin{array}{l}
\lbrack Y_{1},Y_{j}]=Y_{1+j},\quad j=2,...,n-1 \\ 
\lbrack Y_{i},Y_{n-i+1}]=(-1)^{i+1}Y_{n},\quad i=2,...,p
\end{array}
\right. \quad \text{where }n=2p.
\end{equation*}
\end{proposition}

This shows that the preceding algebras, which are also the only naturally graded filiform Lie algebras \cite{Ve}, are of maximal rank. In view of the results in the preceding section we obtain:

\begin{corollary}
Any naturally graded filiform Lie algebra is completable.
\end{corollary}

\begin{theorem1}
Every filiform Lie algebra of rank 1 and dimension n is
isomorphic to one of the following Lie algebras

\begin{enumerate}

\item $A_{n}^{k}\left( \lambda _{1},...,\lambda _{t-1}\right) ,$ $t=\left[ 
\frac{n-k+1}{2}\right] $ , $2\leq k\leq n-3$

\begin{equation*}
\left\{ 
\begin{array}{l}
\left[ Y_{1},Y_{i}\right] =Y_{i+1},\text{ }i=2,...,n-1 \\ 
\left[ Y_{i},Y_{i+1}\right] =\lambda _{i-1}Y_{2i+k-1}\text{ },\text{ }2\leq
i\leq t \\ 
\left[ Y_{i},Y_{j}\right] =a_{ij}Y_{i+j+k-2}\text{ },\text{ }2\leq i\leq j\text{
, }i+j+k-2\leq n
\end{array}
\right.
\end{equation*}

\item $B_{n}^{k}\left( \lambda _{1},...,\lambda _{t-1}\right) $ $n=2m$ , $t=%
\left[ \frac{n-k}{2}\right] $ , $2\leq k\leq n-3$

\begin{equation*}
\left\{ 
\begin{array}{l}
\left[ Y_{1},Y_{i}\right] =Y_{i+1}\text{ }i=2,...,n-2 \\ 
\left[ Y_{i},Y_{n-i+1}\right] =\left( -1\right) _{n}^{i+1}Y\text{ },\text{ }%
i=2,...,n-1 \\ 
\left[ Y_{i},Y_{i+1}\right] =\lambda _{i-1}Y_{2i+k-1}\text{ },\text{ }i=2,...,t
\\ 
\left[ Y_{i},Y_{j}\right] =a_{ij}Y_{i+j-k-2}\text{ },\text{ }2\leq i,j\leq n-2%
\text{ , }i+j+k-2\leq n-2\text{ , }j\neq i+1
\end{array}
\right.
\end{equation*}

\item $C_{n}^{{}}\left( \lambda _{1},...,\lambda _{t}\right) $ , $n=2m+2$ , $%
t=m-1$

\begin{equation*}
\left\{ 
\begin{array}{l}
\left[ Y_{1},Y_{i}\right] =Y_{i+1}\text{ }i=2,...,n-2 \\ 
\left[ Y_{i},Y_{n-i+1}\right] =\left( -1\right) _{n}^{i-1}Y_{n}\text{ },%
\text{ }i=2,...m+1 \\ 
\left[ Y_{i},Y_{n-i-2k+1}\right] =\left( -1\right) ^{i+1}\lambda _{k}Y_{n}%
\text{ },\text{ }i=2,...,n-2-2k\text{ , }k=1,...,m-1
\end{array}
\right.
\end{equation*}
\end{enumerate}

The non defined brackets are equal to zero. In this theorem, $[x]$ denotes
the integer part of $x$ and $\left( \lambda _{1},...,\lambda _{t}\right) $
are non simultaneously vanishing parameters satisfying polynomial equations
associated to the Jacobi conditions. Moreover, the constants $a_{ij}$
satisfy 
\begin{equation*}
a_{ij}=a_{ij+1}+a_{i+1,j}
\end{equation*}
and $a_{ii+1}=\lambda _{i-1}.$
\end{theorem1}

The semi-simple derivations of these algebras have been studied in \cite{GH} and \cite{AG2} in order to determine their maximal tori.

\begin{theorem}
Let $\frak{n}$ be a filiform Lie algebra of rank one. Then $\frak{n}$ is completable if and only if $\frak{n}\simeq A_{n}^{k}\left( \lambda _{1},...,\lambda _{t-1}\right)$ or $\frak{n}\simeq B_{n}^{k}\left( \lambda _{1},...,\lambda _{t-1}\right) $.
\end{theorem}

\begin{proof}
Suppose that $\frak{n}\simeq A_{n}^{k}\left( \lambda _{1},...,\lambda _{t-1}\right)$. From the structure of the algebra it can be seen without difficulty that the basis chosen in the preceding theorem is a basis of eigenvectors for a diagonalizable derivation \cite{GH}. We thus obtain that a maximal torus of derivations $\frak{t}$ is generated by the derivation $f\in Der\left( A_{n}^{k}\left( \lambda _{1},...,\lambda _{t-1}\right)\right)$ defined by
\begin{eqnarray*}
f\left(Y_{1}\right)=Y_{1},\quad f\left(Y_{j}\right)=\left(k-2+j\right)Y_{j},\quad 2\leq j\leq n
\end{eqnarray*}
The weight system is therefore $\left\{\alpha, (k+j-2)\alpha\right\}_{2\leq j\leq n}$. Thus the semidirect product $\frak{t}\oplus\frak{n}$ decomposes as 
\begin{eqnarray*}
\frak{t}\oplus\frak{n}=\frak{t}\oplus\frak{n}_{\alpha}\oplus\sum_{j=2}^{n}(\frak{n}_{(k+j-2)\alpha})
\end{eqnarray*}
and $\left\{\alpha, k\alpha\right\}$ generates $\frak{t}^{*}$. As the algebra $\frak{t}\oplus\frak{n}$ satisfies the conditions of theorem 3, it is complete of non-maximal rank.\newline If $\frak{n}\simeq B_{n}^{k}\left( \lambda _{1},...,\lambda _{t-1}\right) $, the reasoning is quite the same, where a weight system for $B_{n}^{k}\left( \lambda _{1},...,\lambda _{t-1}\right) $ is given by $\left\{\alpha,k\alpha,..,(k+2p+1)\alpha, (2k+2p+1)\alpha\right\}$. Let us therefore suppose that $\frak{n}$ is not isomorphic to the preceding Lie algebras. By theorem ? it must be isomorphic to $ C_{n}^{{}}\left( \lambda _{1},...,\lambda _{t}\right) $. A maximal torus of derivations is generated by the semi-simple derivation $f$ defined by 
\begin{eqnarray*}
f\left(Y_{1}\right)=0,\quad f\left(Y_{j}\right)=Y_{j}, 2\leq j\leq n-1, \quad f\left(Y_{n}\right)=2Y_{n}
\end{eqnarray*}
Observe that condition (ii) of theorem 3 is not satisfied, since there is a zero weight. Consider the endomorphism $h$ defined by $h\left(Y_{2}\right)=Y_{n-1}$ and zero elsewhere. It is obvious that $h\in Der\left(C_{n}^{{}}\left( \lambda _{1},...,\lambda _{t}\right)\right)$, and as $Y_{n-1}\in im(ad(X))$ if and only if $X=Y_{1}$, the derivation is not inner and $C_{n}^{{}}\left( \lambda _{1},...,\lambda _{t}\right) $ cannot be completable.
\end{proof} 

There remain those filiform Lie algebras of rank zero, i.e., those for which the Lie algebra of derivations is itself nilpotent \cite{AC4}. These algebras are definitively non completable since any nilpotent Lie algebra has an outer derivation. We can resume the preceding results in the following 

\begin{theorem}
A filiform Lie algebra of non-zero rank is completable if and only if it is isomorphic to $L_{n}, Q_{n}, A_{n}^{k}\left( \lambda _{1},...,\lambda _{t-1}\right)$ or $B_{n}^{k}\left( \lambda _{1},...,\lambda _{t-1}\right)$.
\end{theorem}

Let $L^{n}$ be the algebraic variety of complex Lie algebra laws on $\mathbb{%
C}^{n}$ \cite{AG1}. Consider the natural action of the algebraic group $GL\left( n,%
\mathbb{C}\right) $ on $L^{n}$ given by 
\begin{eqnarray*}
GL\left( n,\mathbb{C}\right) \times L^{n} &\rightarrow &L^{n} \\
(f,\mu ) &\rightarrow &f\ast \mu
\end{eqnarray*}
with $(f\ast \mu )(X,Y)=f^{-1}(\mu (f(X),f(Y))$ for all $X,Y\in \mathbb{C}%
^{n}.$ We note by $\mathcal{O(\mu )}$ the orbit of $\mu $.

\begin{definition}
The Lie algebra law $\mu $ (or the complex Lie algebra $\frak{g}$ of law $%
\mu )$ is called rigid if $\mathcal{O(\mu )}$ is a Zariski open set of $%
L^{n}.$
\end{definition}

\begin{corollary}
Any filiform Lie algebra which is isomorphic to the nilradical of a solvable rigid law is completable.
\end{corollary}

\section{Complete Lie algebras with $dim H^{2}\left(\frak{r},\frak{r}\right)\geq n$}

In the previous section we have seen which of the filiform Lie algebras are completable. It has therefore sense to ask what happens for the second adjoint cohomology group $H^{2}\left(\frak{n}\oplus\frak{t},\frak{n}\oplus\frak{t}\right)$, where $\frak{n}\oplus\frak{t}$ is complete. We will prove that the dimension of this group can be as high as wanted. Consider the Lie algebra $A_{k+h+3}^{k}\left( \lambda _{1},...,\lambda _{\lbrack\frac{h}{2}\rbrack +1}\right)$ with $k$ even and $\geq 4$. Suppose moreover that $3\leq h\leq k+3$. Let us consider the basis $\left\{Y_{1},Y_{k},..,Y_{2k+2h+1}\right\}$, where the subindex makes reference to the weight of the vector by the action of the torus $\frak{t}$. The brackets of the algebra are given by
\begin{eqnarray*}
\lbrack Y_{1},Y_{k+i}\rbrack = Y_{k+1+i}, 0\geq i\geq 2k+2h\\
\lbrack X_{k+i},X_{k+j}\rbrack=a_{ij}X_{2k+i+j}, 0\geq i,j\geq \lbrack\frac{h}{2}\rbrack+1
\end{eqnarray*}
We have that the coefficients $a_{ij}$ satisfy the equations

\begin{eqnarray}
a_{ij}=a_{i,j-1}-a_{i+1,j-1}\\
a_{i,i+1}=\lambda_{i+1}
\end{eqnarray}
for $0\geq i<j\geq \lbrack\frac{p}{2}\rbrack$.\newline The preceding relations are an immediate consequence of the Jacobi conditions related to the triples $\left\{X_{1},X_{i},X_{j}\right\}$, which are the only ones since $3k+3\geq 2k-2+p$. It is not difficult to see that, whenever $\lambda_{1}$ is non-zero, it can be normalized. In this manner we obtain $\lbrack\frac{h}{2}\rbrack$ essential parameters which separate isomorphism classes. From equations $(1)$ and $(2)$ we have that 
\begin{eqnarray*}
a_{ij}=\sum_{k=1}^{\lbrack\frac{p}{2}\rbrack+1}\alpha_{ij}^{k}\lambda_{k}
\end{eqnarray*}
where the $\alpha_{ij}^{k}$ are the integers coefficients obtained by the resolution of the systems associated to the Jacobi conditions giving the equations. If we define the functions 
\begin{eqnarray*}
h_{k}\left(a_{ij}\right)=\alpha_{ij}^{k}, 
\end{eqnarray*}
we see that the bilinear alternated mapping $\varphi_{k}$ defined by 
\begin{eqnarray*}
\left( X_{k+i},X_{k+j}\right)\rightarrow h_{k}\left(a_{ij}\right)X_{2k+i+j}
\end{eqnarray*}
defines a linear deformation \cite{G}. Moreover, the $\varphi_{k}$ cannot be coboundaries since otherwise the parameters $\lambda_{k}$ would not separate isomorphism classes. We obtain $\varphi_{k}\in Z^{2}\left(\frak{r}_{h},\frak{r}_{h}\right)$, where $\frak{r}_{h}=A_{k+3+h}^{k}\left( \lambda _{1},...,\lambda _{t-1}\right)\oplus\frak{t}$. We obtain 

\begin{theorem}
In the preceding conditions, the complete Lie algebra  $\frak{r}_{h}=A_{k+3+h}^{k}\left( \lambda _{1},...,\lambda _{\lbrack\frac{h}{2}\rbrack+1}\right)\oplus\frak{t}$ satisfies $dim H^{2}\left(\frak{r}_{h},\frak{r}_{h}\right)\geq \lbrack\frac{h}{2}\rbrack$.
\end{theorem}

\begin{corollary}
For any $n\in\mathbb{N}^{+}$ there exists a complete Lie algebra $\frak{r}$ such that $dim H^{2}\left(\frak{r},\frak{r}\right)\geq n$.
\end{corollary}  

The reasoning just described is essentially the same applied by Carles in \cite{C} to construct the first examples of solvable rigid laws with non-vanishing second cohomology group. We remark that the examples presented by him in \cite{C} were also complete Lie algebras.

\end{document}